\documentclass[12pt,a4paper,reqno]{amsart}
\usepackage{graphics,epic}
\usepackage{amsmath,amssymb, amsthm}
\usepackage[all,2cell]{xy}

\newtheorem{theorem}{Theorem}[section]
\newtheorem*{theorem*}{Theorem}

\newtheorem{lemma}[theorem]{Lemma}
\newtheorem{proposition}[theorem]{Proposition}
\newtheorem{corollary}[theorem]{Corollary}

\newtheorem*{conjecture*}{Conjecture}

\newtheorem{remark}[theorem]{Remark}

\renewcommand{\hat}[1]{\widehat{#1}}

\newcommand{\id}{{\rm id}}
\newcommand{\ch}{{\rm ch}}
\newcommand{\im}{{\rm im}}
\newcommand{\Hom}{{\rm Hom}\,}
\newcommand{\Glob}{{\rm Glob}\,}
\newcommand{\End}{{\rm End}\,}

\newcommand{\qdim}{{\rm qdim}\,}

\newcommand{\Z}{\mathbb{Z}}

\newcommand{\C}{\mathbb{C}}

\def\wt{{\rm wt}}
\def\C{{\mathbb C}}

\def\Z{{\mathbb Z}}

\def\1{{\bf 1}}
\def\tr{{\rm tr}}
\def \End{{\rm End}}
\def \Hom{{\rm Hom}}

\def \Ind{{\rm Ind}}

\def \pf{\noindent {\bf Proof: \,}}
\def\theequation{5.\arabic{equation}}
\def \h{\mathfrak{h}}
\def \w{\omega}
\def \g{\mathfrak{g}}
\def \p{\mathfrak{p}}

\begin{document}

\title[Quantum dimensions and irreducible modules of  diagonal coset VOAs]{Quantum dimensions and irreducible modules of some diagonal coset vertex operator algebras}

\author{Xingjun Lin}
\address{Xingjun Lin, Collaborative innovation centre of Mathematics, School of Mathematics and Statistics, Wuhan University, Luojiashan, Wuhan, Hubei, 430072, China.}
\thanks{X. Lin was supported by China NSF grant
11801419 and the starting research fund from Wuhan University (No. 413100051)}
\email{linxingjun88@126.com}
\begin{abstract}
In this paper,  under the assumption that the diagonal coset vertex operator algebra $C(L_{\g}(k+l,0),L_{\g}(k,0)\otimes L_{\g}(l,0))$ is rational and $C_2$-cofinite, the global dimension of $C(L_{\g}(k+l,0),L_{\g}(k,0)\otimes L_{\g}(l,0))$ is obtained,  the quantum dimensions of multiplicity spaces viewed as $C(L_{\g}(k+l,0),L_{\g}(k,0)\otimes L_{\g}(l,0))$-modules are also obtained. As an application, a method to classify irreducible modules of $C(L_{\g}(k+l,0),L_{\g}(k,0)\otimes L_{\g}(l,0))$ is provided. As an example,  we prove that the diagonal coset vertex operator algebra $C(L_{E_8}(k+2,0),L_{E_8}(k,0)\otimes L_{E_8}(2,0))$ is  rational, $C_2$-cofinite, and classify irreducible modules of $C(L_{E_8}(k+2,0),L_{E_8}(k,0)\otimes L_{E_8}(2,0))$.
\end{abstract}
\maketitle
\section{Introduction \label{intro}}
\def\theequation{1.\arabic{equation}}
\setcounter{equation}{0}

Let $\g$ be a reductive finite dimensional Lie algebra and $\p$ be a reductive subalgebra of $\g$. It is important to study  the multiplicity spaces of an integrable highest weight representation of the affine Lie algebra $\hat \g$ with respect to the affine subalgebra $\hat \p$. The conformal invariance properties of the multiplicity spaces of an integrable highest weight representation of $\hat \g$ with respect to $\hat \p$ were first studied in \cite{GKO}, unitary highest weight modules of the Virasoro algebra with central charge less than 1 were constructed in this important paper.  Later, the modular invariance properties of the branching functions  of an integrable highest weight representation of $\hat \g$ with respect to $\hat \p$ were studied in important papers  \cite{KP}, \cite{KW}.

For a positive integer $k$, it is well-known that the vacuum module $L_{\g}(k, 0)$ of $\hat \g$ is a vertex operator algebra \cite{FZ}. Since $\hat \p$ is a subalgebra of $\hat \g$, $L_{\g}(k, 0)$ may be viewed as a $\hat \p$-module of level $l$, where $l$ is a positive integer. It was realized in \cite{FZ} that the multiplicity space $C(L_{\p}(l, 0), L_{\g}(k, 0)) $ of $L_{\p}(l, 0)$ in $L_{\g}(k, 0)$ is also a vertex operator algebra. Moreover, the multiplicity spaces of an integrable highest weight representation of $\hat \g$ with respect to $\hat \p$ are modules of the vertex operator algebra $C(L_{\p}(l, 0), L_{\g}(k, 0)) $. In general, the multiplicity spaces of an integrable highest weight representation of $\hat \g$ with respect to $\hat \p$ may be reducible $C(L_{\p}(l, 0), L_{\g}(k, 0))$-modules. In addition, viewed as $C(L_{\p}(l, 0), L_{\g}(k, 0))$-modules, there may be identifications between  the multiplicity spaces of integrable highest weight representations of $\hat \g$ with respect to $\hat \p$. Thus, it is interesting to study the structures of the multiplicity spaces as $C(L_{\p}(l, 0), L_{\g}(k, 0))$-modules.

There is another approach to study the multiplicity spaces of an integrable highest weight representation of $\hat \g$ with respect to $\hat \p$ in the framework of conformal nets \cite{X}. Many important results have been established in \cite{X}. However, it seems that the method in \cite{X} cannot be used to study the multiplicity spaces of an admissible highest weight representation of $\hat \g$ with respect to $\hat \p$ \cite{KW2}, because nonunitary representations will appear in this situation. So it is still interesting to study the problem in the framework of vertex operator algebras.

When $\g$ is a simple Lie algebra and $\p$ is a Cartan subalgebra of $\g$. Then the vertex operator algebra $C(L_{\p}(k, 0), L_{\g}(k, 0))$ is the parafermion vertex operator algebra \cite{DLY}. The parafermion vertex operator algebras have been well-studied (see \cite{DKR} and reference therein). In this paper, by considering $\g$ as the diagonal subalgebra of $\g\oplus \g$, we will study the vertex operator algebra $C(L_{\g}(k+l, 0), L_{\g}(k,0)\otimes L_{\g}(l,0))$, which can realize many important  vertex operator algebras (see \cite{ACL}).

First, under the assumption that the vertex operator algebra $C(L_{\g}(k+l,0),L_{\g}(k,0)\otimes L_{\g}(l,0))$ is rational and $C_2$-cofinite, the global dimension of $C(L_{\g}(k+l,0),L_{\g}(k,0)\otimes L_{\g}(l,0))$ is obtained,  the quantum dimensions of the multiplicity spaces viewed as $C(L_{\g}(k+l,0),L_{\g}(k,0)\otimes L_{\g}(l,0))$-modules are also obtained (see Theorem \ref{qdimc}). These results are further used to classify irreducible modules of $C(L_{\g}(k+l,0),L_{\g}(k,0)\otimes L_{\g}(l,0))$ (see Theorem \ref{general}). The method in Theorem \ref{general} to classify irreducible modules of $C(L_{\g}(k+l,0),L_{\g}(k,0)\otimes L_{\g}(l,0))$ depends on the assumption that $C(L_{\g}(k+l,0),L_{\g}(k,0)\otimes L_{\g}(l,0))$ is rational and $C_2$-cofinite. There is another work \cite{CWX} in which $C(L_{\g}(k+l,0),L_{\g}(k,0)\otimes L_{\g}(l,0))$-modules are studied without the assumption.

As an example, we show that the vertex operator algebra $C(L_{E_8}(k+2,0),L_{E_8}(k,0)\otimes L_{E_8}(2,0))$ is rational and $C_2$-cofinite (see Theorem \ref{main1}). We then use the method in Theorem \ref{general} to show that the multiplicity spaces are irreducible modules of $C(L_{E_8}(k+2,0),L_{E_8}(k,0)\otimes L_{E_8}(2,0))$ (see Theorem \ref{main2}). Moreover, we show that the multiplicity spaces are all the non-isomorphic irreducible modules of $C(L_{E_8}(k+2,0),L_{E_8}(k,0)\otimes L_{E_8}(2,0))$ (see Theorem \ref{main2}).

The paper is organized as follows: In Section 2, we recall some basic facts about vertex operator algebras. In Section 3, we prove the diagonal coset vertex operator algebra $C(L_{E_8}(k+2,0),L_{E_8}(k,0)\otimes L_{E_8}(2,0))$ is  rational, $C_2$-cofinite. In Section 4, under the assumption that the diagonal coset vertex operator algebra $C(L_{\g}(k+l,0),L_{\g}(k,0)\otimes L_{\g}(l,0))$ is rational and $C_2$-cofinite, the global dimension of $C(L_{\g}(k+l,0),L_{\g}(k,0)\otimes L_{\g}(l,0))$ is obtained. The quantum dimensions of multiplicity spaces viewed as $C(L_{\g}(k+l,0),L_{\g}(k,0)\otimes L_{\g}(l,0))$-modules are also obtained. In Section 5, a method to classify irreducible modules of $C(L_{\g}(k+l,0),L_{\g}(k,0)\otimes L_{\g}(l,0))$ is provided. As an example,  irreducible modules of $C(L_{E_8}(k+2,0),L_{E_8}(k,0)\otimes L_{E_8}(2,0))$ are classified.
\section{Preliminaries }
\def\theequation{2.\arabic{equation}}
\setcounter{equation}{0}
\subsection{Basics}
In this subsection we briefly review  some basic notions and facts in the theory of vertex operator algebras from \cite{DLM}, \cite{FHL}, \cite{FLM},  \cite{LL} and \cite{Z}. Let $(V, Y, \1, \w)$ be a vertex operator algebra as defined in \cite{FLM} (see also \cite{B}). $V$ is called {\em $C_2$-cofinite} if $\dim V/C_2(V)<\infty$, where $C_2(V)=\langle u_{-2}v|u, v\in V \rangle$. $V$ is called of {\em CFT} type if $V$ has the decomposition $V=\bigoplus_{n\geq 0}V_n$ with respect to $L(0)$ such that $\dim V_0=1$. If $v\in V_n$, the {\em conformal weight} $\wt v$ of $v$ is defined to be $n$.

A {\em weak}  $V$-module $M$ is a vector space equipped
with a linear map
\begin{align*}
Y_{M}:V&\to (\End M)[[z, z^{-1}]],\\
v&\mapsto Y_{M}(v,z)=\sum_{n\in\Z}v_nz^{-n-1},\,v_n\in \End M
\end{align*}
satisfying the following conditions: For any $u\in V,\ v\in V,\ w\in M$ and $n\in \Z$,
\begin{align*}
&\ \ \ \ \ \ \ \ \ \ \ \ \ \ \ \ \ \ \ \ \ \ \ \ \ \ u_nw=0 \text{ for } n>>0;\\
&\ \ \ \ \ \ \ \ \ \ \ \ \ \ \ \ \ \ \ \ \ \ \ \ \ \ Y_M(\1, z)=\id_M;\\
\begin{split}
&z_{0}^{-1}\delta\left(\frac{z_{1}-z_{2}}{z_{0}}\right)Y_{M}(u,z_{1})Y_M(v,z_{2})-z_{0}^{-1}\delta\left(
\frac{z_{2}-z_{1}}{-z_{0}}\right)Y_M(v,z_{2})Y_M(u,z_{1})\\
&\quad=z_{2}^{-1}\delta\left(\frac{z_{1}-z_{0}}{z_{2}}\right)Y_M(Y(u,z_{0})v,z_{2}).
\end{split}
\end{align*}

A weak
 $V$-module  $M$ is called an \textit{admissible} $V$-module if $M$ has a $\Z_{\geq
0}$-gradation $M=\bigoplus_{n\in\Z_{\geq 0}}M(n)$ such
that
\begin{align*}\label{AD1}
a_mM(n)\subset M(\wt{a}+n-m-1)
\end{align*}
for any homogeneous $a\in V$ and $m,\,n\in\Z$.

 A vertex operator algebra $V$ is called {\em simple} if $V$ viewed as a weak $V$-module is irreducible. And $V$ is said to be \textit{rational} if
any  admissible $V$-module is completely reducible. It was proved in \cite{DLM1} that if $V$ is rational then there are only finitely many irreducible admissible $V$-modules up to isomorphism.

An ({\em ordinary})  $V$-module is a weak $V$-module $M$ such that $M=\bigoplus_{\lambda\in\C}
M_{\lambda}$, where
$M_\lambda=\{w\in M|L(0)w=\lambda w\}$, and that $M_\lambda$ is
finite dimensional and for fixed $\lambda\in\C$, $M_{\lambda+n}=0$
for sufficiently small integer $n$.
Let $M = \bigoplus_{\lambda\in \mathbb{C}}{M_{\lambda}}$ be a $V$-module. Set $M'
= \bigoplus_{\lambda \in \mathbb{C}}{M_\lambda^*}$, the restricted
dual of $M$. It was proved in \cite{FHL} that $M'$ is also a
$V$-module such that the vertex operator map $Y'$ is defined
by
$$\langle Y'(a, z)u', v\rangle  = \langle u', Y(e^{zL(1)}(-z^{-2})^{L(0)}a, z^{-1})v\rangle ,$$for $a\in V, u'\in
M'$ and $v\in M$. The $V$-module $M'$ is called the {\em contragredient
module} of $M$. It was also proved in \cite{FHL} that if $M$ is irreducible,  then so
is $M'$, and that $(M')'\simeq M$. A $V$-module $M$ is called {\em self-dual} if $M\cong M'$.
A vertex operator algebra $V$ is called {\em strongly regular} if $V$ is of CFT type, simple, self-dual, rational and $C_{2}$-cofinite.
\subsection{Modular invariance of trace functions of vertex operator algebras}
We now turn our discussion to the modular invariance property in the theory of vertex operator algebras.  Let $\mathcal{H} =\{\tau\in \mathbb{C}| \im\tau>0\}$. Recall that the full modular group $SL(2, \mathbb{Z})$ has generators $S=\left(\begin{array}{cc}0 & -1\\ 1 & 0\end{array}\right)$, $T=\left(\begin{array}{cc}1 & 1\\ 0 & 1\end{array}\right)$ and acts on $\mathcal H$ as follows:$$\gamma: \tau\longmapsto \frac{a\tau+b}{c\tau+d}, \  \gamma=\left(\begin{tabular}{cc}
$a$ $b$\\
$c$ $d$\\
\end{tabular}\right) \in SL(2, \mathbb{Z}).$$

Let $V$ be a rational and $C_2$-cofinite  vertex operator algebra, $M^0, M^1,...,M^p$ be all the irreducible $V$-modules. Then $M^i, 0\leq i \leq p$, has the form $$M^i=\bigoplus_{n=0}^{\infty}M^i_{\lambda_i+n},$$ with $M^i_{\lambda_i}\neq 0$ for some number $\lambda_i$ which is called {\em conformal weight} of $M^i$.
  For any irreducible $V$-module $M^i$, the {\em trace function} associated to $M^i$ is defined as follows: For any homogenous element $v\in V$ and $\tau\in \mathcal H$,
\begin{equation*}
Z_{M^i}(v,\tau):=\tr_{M^i}o(v)q^{L(0)-c/24}=q^{\lambda_i-c/24}\sum_{n\in\mathbb{Z}_{\geq 0}} \tr_{M^i_{\lambda_i+n}}o(v)q^n,
\end{equation*}
where $o(v)=v(\wt v-1)$ and $q=e^{2\pi \sqrt{-1}\tau}$. Since $V$ is $C_2$-cofinite, $Z_{M^i}(v,\tau)$ converges to a holomorphic function on the domain $|q| < 1$ \cite{DLM2}, \cite{Z}. The $Z_M(\1, \tau)$ which is also denoted by $\ch_q M$ is called the $q$-character of $M$.

Another vertex operator algebra structure $(V, Y[., z], \1, \w-c/24)$ is defined on $V$ in \cite{Z} with grading
$$V=\oplus_{n\geq 0} V_{[n]}.$$
For $v\in V_{[n]}$, we write $\wt[v]=n$.
\begin{theorem}\label{minvariance}
 Let $V$ be a rational and $C_2$-cofinite vertex operator algebra with the irreducible $V$-modules $M^0,...,M^p.$   Then there is a group homomorphism $\rho_V: SL(2, \Z)\to GL_{p+1}(\C)$ with $\rho(\gamma)=(\gamma_{i,j})$ such that for any $0\leq i\leq p$ and homogeneous $v\in V_{[n]}$, $$Z_{M^i}(v,\frac{a\tau+b}{c\tau+d})=(c\tau+d)^n\sum_{j=0}^p \gamma_{i,j}Z_{M^j}(v,\tau).$$
 Moreover, the matrix $(\gamma_{i,j})$ is independent of the choice of $v\in V$.
\end{theorem}
The theorem was proved in \cite{Z} (also see \cite{DLM2}). We will also use $S$ to denote the  matrix $\rho_V(S)=(S_{i,j})$. It is known that $S$ is a symmetric matrix.
\subsection{Quantum dimensions}
In this subsection,  we recall some facts about quantum dimensions of modules of vertex operator algebras from \cite{DJX}. Let $V$ be a strongly regular vertex operator algebra and let $M^0=V, M^1,...,M^p$ be all the inequivalent irreducible $V$-modules. For a $V$-module $M$, the {\em quantum dimension} of $M$ is defined to be
\begin{align*}
\qdim_V M=\lim_{y\to 0^+}\frac{Z_{M}(\1, \sqrt{-1}y)}{Z_{V}(\1, \sqrt{-1}y)},
\end{align*}
where $y$ is real and positive. The {\em global dimension}  of the vertex operator algebra $V$ is defined to be
\begin{align*}
{\rm Glob}\, V=\sum_{i=0}^p(\qdim_V M^i)^2.
\end{align*}

The following result was proved in  \cite{DJX}.
\begin{theorem}\label{qdim1}
Let $V$ be a strongly regular vertex operator algebra and let $M^0=V, M^1,...,M^p$ be all the irreducible $V$-modules. Assume further that the conformal weights of $ M^1,...,M^p$ are greater than $0$. Then \\
(1) $\qdim M^i\geq 1$ for any $0\leq i\leq p$.\\
(2) $\qdim_V M^i=\frac{S_{0,i}}{S_{0,0}}$.\\
(3) ${\rm Glob}\, V=\frac{1}{S_{0,0}^2}.$\\
(4) $M^i$ is a simple current $V$-module if and only if $\qdim M^i=1$.
\end{theorem}
\section{Rationality of some coset vertex operator algebras}
\def\theequation{3.\arabic{equation}}
\setcounter{equation}{0}
 \subsection{Affine vertex operator algebras}  In this subsection, we shall recall some facts about affine vertex operator algebras from \cite{FZ} and \cite{LL}. Let $\g$ be a finite dimensional simple Lie algebra and $\langle\, ,\, \rangle$ the normalized Killing form of $\g$, i.e., $\langle\theta, \theta\rangle=2$ for the highest root $\theta$ of $\g$. Fix a Cartan subalgebra $\h$ of $\g$ and  denote the corresponding root system by $\Delta_{\g}$ and the root lattice by $Q$. We further fix simple roots $\{\alpha_1,\cdots,\alpha_l\}$, and denote the set of positive roots by $\Delta_{\g}^+$. Then the weight lattice $P$ of $\g$ is the set of $\lambda\in \h$ such that $\frac{2\langle\lambda, \alpha\rangle}{\langle\alpha, \alpha\rangle}\in\Z$ for all $\alpha\in \Delta_{\g}$. Note that $P$ is equal to $\oplus_{i=1}^l\Z\Lambda_i$, where $\Lambda_i$ are the fundamental weights defined by the equation $\frac{2\langle\Lambda_i, \alpha_j\rangle}{\langle\alpha_j, \alpha_j\rangle}=\delta_{i,j}$. We also use the standard notation $P_+$ to denote the set of dominant weights $\{\Lambda\in P\mid\frac{2\langle\Lambda, \alpha_j\rangle}{\langle\alpha_j, \alpha_j\rangle}\geq 0,~1\leq j\leq l \}$. For any  $\alpha\in \Delta_{\g}^+$, we fix $x_{\pm\alpha}\in \g_{\pm\alpha}$ such that $[x_{\alpha},x_{-\alpha}]=h_{\alpha}$, $[h_{\alpha},x_{\pm\alpha}]=\pm 2 x_{\pm\alpha}$, where $h_{\alpha}=\frac{2}{\langle\alpha,\alpha\rangle}\alpha$.

Recall that the affine Lie algebra associated to $\g$ is defined on $\hat{\g}=\g\otimes \C[t^{-1}, t]\oplus \C K$ with Lie brackets
\begin{align*}
[x(m), y(n)]&=[x, y](m+n)+\langle x, y\rangle m\delta_{m+n,0}K,\\
[K, \hat\g]&=0,
\end{align*}
for $x, y\in \g$ and $m,n \in \Z$, where $x(n)$ denotes $x\otimes t^n$.

For a positive integer $k$ and a weight $\Lambda \in P$, let $L_{\g}(\Lambda)$ be the irreducible highest weight module for $\g$ with highest weight $\Lambda$ and define
\begin{align*}
V_{\g}(k, \Lambda)=\Ind_{\g\otimes \C[t]\oplus \C K}^{\hat \g}L_{\g}(\Lambda),
\end{align*}
where $L_{\g}(\Lambda)$ is viewed as a module for $\g\otimes \C[t]\oplus \C K$ such that $\g\otimes t\C[t]$ acts as $0$ and $K$ acts as $k$. It is well-known that $V_{\g}(k, \Lambda)$ has a unique maximal proper submodule which is denoted by $J(k, \Lambda)$ (see \cite{K}). Let $L_{\g}(k, \Lambda)$ be the corresponding irreducible quotient module. It was proved in \cite{FZ} that $L_{\g}(k, 0)$ has a vertex operator algebra structure such that the Virasoro vector
\begin{align*}\label{virasoro}
\w=\frac{1}{2(k+h^{\vee})} \sum_{i=1}^{\dim \g} u_i(-1)u_i(-1)\1,
\end{align*}
where $h^\vee$ denotes the dual Coxeter number of $\g$ and $\{u_i|1\leq i\leq \dim \g\}$ is an orthonormal basis of $\g$ with respect to $\langle,\rangle$.
Moreover, we have the following results which were proved in \cite{DLM2}, \cite{FZ}, \cite{K}.
\begin{theorem}
Let $k$ be a positive integer. Then \\
(1) $L_{\g}(k, 0)$ is a strongly regular vertex operator algebra.\\
 (2) $L_{\g}(k, \Lambda)$ is a module for the vertex operator algebra $L_{\g}(k, 0)$ if and only if $\Lambda \in P_+^k$, where $P_+^k=\{\Lambda \in P_+|\langle\Lambda, \theta\rangle\leq k\}$.\\
 (3) If $L_{\g}(k, \Lambda)$ is an $L_{\g}(k, 0)$-module  such that $L_{\g}(k, \Lambda)\ncong L_{\g}(k, 0)$, then the conformal weight of $L_{\g}(k, \Lambda)$ is positive.
\end{theorem}

We next recall some fact about modular invariance properties about $L_{\g}(k, 0)$. Since $L_{\g}(k, 0)$ is strongly regular and $\{L_{\g}(k, \Lambda)|\Lambda \in P_+^k\}$ are all the non-isomorphic irreducible $L_{\g}(k, 0)$-modules, the trace functions  $\{Z_{L_{\g}(k, \Lambda)}(v,\tau)|\Lambda \in P_+^k\}$ is closed under the action of $SL(2,\Z)$.
\begin{theorem}\cite{K,KW}
 For $\Lambda, \Lambda' \in P_+^k$, let $S_{L_{\g}(k, \Lambda), L_{\g}(k, \Lambda')}$ be complex numbers such that
 $$Z_{L_{\g}(k, \Lambda)}(v,\frac{-1}{\tau})=\tau^{\wt [v]}\sum_{\Lambda' \in P_+^k}S_{L_{\g}(k, \Lambda), L_{\g}(k, \Lambda')}Z_{L_{\g}(k, \Lambda')}(v,\tau).$$
 Then
$$S_{L_{\g}(k, 0), L_{\g}(k, \Lambda)}=|P/(k+h^\vee)Q_L|^{-1/2} (k+h^\vee)^{-l/2}\prod_{\alpha\in \Delta_{\g}^+}2\sin \frac{\pi \langle\Lambda +\rho, \alpha\rangle}{k+h^\vee},$$
where $Q_L\subseteqq Q$ is the sublattice of $Q$ spanned by all long roots and $\rho=\sum_{i=1}^l\Lambda_i$.
\end{theorem}

 Finally, we recall some facts about simple current modules of $L_{\g}(k, 0)$. Let $\theta=\sum_{i=1}^l a_i\alpha_i$, $a_i\in \Z_+$, be the highest root. It is well-known that the irreducible $L_{\g}(k, 0)$-module $L_{\g}(k, k\Lambda_i)$ is a simple current $L_{\g}(k, 0)$-module if $a_i=1$ (see \cite{DLM0,F,FG,L3}). Set $J=\{i|a_i=1\}\cup\{0\}$. Then it is known that $|J|=|P/Q|$ \cite{KW}. In particular,  for $\g=E_8$, $J=\{0\}$ \cite{KW}.
\subsection{Virasoro vertex operator algebras} In this subsection, we recall some facts about
Virasoro vertex operator algebras \cite{FZ}, \cite{W}.  Let $L=\oplus_{n\in \mathbb{Z}}
\mathbb{C}L_n\oplus \mathbb{C}C$ be the Virasoro algebra with the
commutation relations
\begin{align*}
[L_m, L_n]=(m-n)&L_{m+n}+\frac{1}{12}(m^3-m)\delta_{m+n,
0}C,\\
&[L_m, C]=0.
\end{align*}

Consider the subalgebra $\mathfrak{b}=(\oplus_{n\geq
1}\mathbb{C}L_n)\oplus(\mathbb{C}L_0\oplus \mathbb{C}C)$ of $L$. For any two
complex numbers $c, h\in \C$, let $\C$ be a 1-dimensional
$\mathfrak{b}$-module such that:
$$L_n\cdot 1=0, n\geq 1,~
L_0\cdot 1=h\cdot 1,~C\cdot 1=c\cdot 1.
$$
Define $V(c, h)=U(L)\otimes_{U(\mathfrak{b})}\C$, where
$U(\cdot)$ denotes the universal enveloping algebra. Then $V(c, h)$ is
a highest weight module of the Virasoro algebra of highest weight $(c, h)$,
and $V(c, h)$ has a unique maximal proper submodule $J(c, h)$. Let
$L(c, h)$ be the unique irreducible quotient module of $V(c, h)$.
Set
$\overline{V(c, 0)}=V(c, 0)/(U(L)L_{-1}1\otimes 1),$ it is well-known
that $\overline{V(c, 0)}$ has a vertex operator algebra structure with Virasoro vector $\w=L_{-2}1$
and $ L(c, 0)$ is the unique irreducible quotient vertex operator
algebra of $\overline{V(c, 0)}$ \cite{FZ}.

For coprime integers $p, q\geq 2$, set
\begin{align*}
c_{p,q}=1-\frac{6(p-q)^2}{pq},~~~
h^{p,q}_{r,
s}=\frac{(rp-sq)^2-(p-q)^2}{4pq},~ (r,s)\in E_{p,q},
 \end{align*}
 where   $E_{p,q}$ denotes the set $\{(r, s)|1\leq s\leq p-1,~1\leq r \leq q-1,~r+s\equiv 0~{\rm mod}~2\}$.
\begin{theorem}\cite{DLM2,W}
 For coprime integers $p, q\geq 2$, the Virasoro vertex operator algebra $L(c_{p,q}, 0)$ is strongly regular and $L(c_{p,q}, h^{p,q}_{r, s}),~ (r,s)\in E_{p,q},$ are
the complete list of irreducible $L(c_{p,q}, 0)$-modules.
\end{theorem}

In particular, when $p=3, q=4$, the Virasoro vertex operator algebra $L(1/2, 0)$ has three irreducible modules $L(1/2, 0), L(1/2, 1/2), L(1/2, 1/16)$. In the following, we will need the fusion rules \cite{FHL} between irreducible modules of $L(1/2, 0)$.
\begin{theorem}\cite{DMZ, W}\label{fusion}
The fusion rules between irreducible $L(1/2, 0)$-modules are as follows:
\begin{align*}
&L(1/2, 1/2)\times L(1/2, 1/2)=L(1/2, 0),\\
&L(1/2, 1/16)\times L(1/2, 1/2)=L(1/2, 1/16),\\
&L(1/2, 1/16)\times L(1/2, 1/16)=L(1/2, 0)\oplus L(1/2, 1/2).
\end{align*}
\end{theorem}
\subsection{Rationality of coset vertex operator algebra $C(L_{E_8}(k+2,0), L_{E_8}(k,0)\otimes L_{E_8}(2,0))$} Let $\g$ be a finite dimensional simple Lie algebra, $k,l$ be positive integers. Consider the vertex operator algebra $L_{\g}(k,0)\otimes L_{\g}(l,0)$, and the vertex subalgebra $W$ of  $L_{\g}(k,0)\otimes L_{\g}(l,0)$ generated by $\{x(-1)\1\otimes \1+\1\otimes x(-1)\1|x\in \g\}$. By Theorem 3.1 of \cite{DM2}, $W$ is isomorphic to $L_{\g}(k+l,0)$. Set \begin{align*}
C(L_{\g}&(k+l,0),L_{\g}(k,0)\otimes L_{\g}(l,0))
\\&=\{u\in L_{\g}(k,0)\otimes L_{\g}(l,0)|u_{n}v=0, \forall v\in L_{\g}(k+l,0), \forall n\in \Z_{\geq 0}\}.
\end{align*}
It is well-known that $C(L_{\g}(k+l,0),L_{\g}(k,0)\otimes L_{\g}(l,0))$ is a vertex subalgebra of $L_{\g}(k,0)\otimes L_{\g}(l,0)$ \cite{LL}. Moreover, by
Proposition 12.10 of \cite{K} and Lemma 2.1 of \cite{ACKL}, we have
\begin{lemma}\label{scoset}
Let $\w^1, \w^2, \w^3$ be the Virasoro elements of $L_{\g}(k,0)$, $L_{\g}(l,0)$, $L_{\g}(k+l,0)$, respectively. Then $C(L_{\g}(k+l,0),L_{\g}(k,0)\otimes L_{\g}(l,0))$  is a  simple vertex operator algebra with the Virasoro vector $\w^1+\w^2-\w^3$.
\end{lemma}
In the following, we consider the case $\g$ is the simple Lie algebra of type $E_8$. We shall prove that the coset vertex operator algebra $C(L_{E_8}(k+2,0), L_{E_8}(k,0)\otimes L_{E_8}(2,0))$ is rational and $C_2$-cofinite. To prove this result, we need the following important result which was obtained in \cite{ACL}.
\begin{theorem}\label{coset1}
Let $k$ be a positive integer. Then the coset vertex operator algebra $C(L_{E_8}(k+1,0), L_{E_8}(k,0)\otimes L_{E_8}(1,0))$ is strongly regular.
\end{theorem}

We also need the following facts were obtained in the formulas (4.1.8), (4.1.9) of \cite{KW}.
\begin{lemma}\label{coset2}
(1) $C(L_{E_8}(2, 0), L_{E_8}(1, 0)\otimes L_{E_8}(1, 0))$ is isomorphic to the Virasoro vertex operator algebra $L(\frac{1}{2}, 0)$.
\\(2) $L_{E_8}(1, 0)\otimes L_{E_8}(1, 0)$ viewed as an $L_{E_8}(2, 0)\otimes L(\frac{1}{2}, 0)$-module has the following decomposition
$$L_{E_8}(1, 0)\otimes L_{E_8}(1, 0)\cong L_{E_8}(2, 0)\otimes L(\frac{1}{2}, 0)\oplus L_{E_8}(2, \Lambda_7)\otimes L(\frac{1}{2}, \frac{1}{2}) \oplus L_{E_8}(2, \Lambda_1)\otimes L(\frac{1}{2}, \frac{1}{16}).$$
\end{lemma}

\vskip .25cm
 Recall that a vertex operator algebra $U$ is called an {\em extension} of $V$ if $V$ is a vertex operator subalgebra of $U$ and $V$, $U$ have the same Virasoro vector. We will also need the following results about extensions.
    \begin{theorem}\label{abd}\cite{ABD}
Let $V$ be a $C_2$-cofinite vertex operator algebra and $U$ be an extension of $V$. Assume further that $U$ viewed as a $V$-module is completely reducible, then $U$ is $C_2$-cofinite.
\end{theorem}

  \begin{theorem}\label{cvoa2}  \cite {HKL}
Let $V$ be a strongly regular vertex operator algebra. Suppose that $U$ is a simple vertex operator algebra and an extension of $V$, then $U$ is rational.
\end{theorem}

We are now ready to prove the main result in this subsection.
\begin{theorem}\label{main1}
Let $k$ be a positive integer. Then the coset vertex operator algebra $C(L_{E_8}(k+2, 0), L_{E_8}(k, 0)\otimes L_{E_8}(2, 0))$ is rational and $C_2$-cofinite.
\end{theorem}
\pf Consider the vertex operator algebra $L_{E_8}(k, 0)\otimes L_{E_8}(1, 0)\otimes L_{E_8}(1, 0)$. By the diagonal embedding, $L_{E_8}(k+2, 0)$ can be viewed as a subalgebra of $L_{E_8}(k, 0)\otimes L_{E_8}(1, 0)\otimes L_{E_8}(1, 0)$. Consider the coset vertex operator algebra $C(L_{E_8}(k+2, 0),L_{E_8}(k, 0)\otimes L_{E_8}(1, 0)\otimes L_{E_8}(1, 0))$. It is clear that $C(L_{E_8}(k+2, 0),L_{E_8}(k, 0)\otimes L_{E_8}(1, 0)\otimes L_{E_8}(1, 0))$ is an extension of $$C(L_{E_8}(k+1, 0),L_{E_8}(k, 0)\otimes L_{E_8}(1, 0))\bigotimes C(L_{E_8}(k+2, 0),L_{E_8}(k+1, 0)\otimes L_{E_8}(1, 0)).$$
Thus, by Lemma \ref{scoset} and Theorems \ref{coset1}, \ref{abd}, \ref{cvoa2},  $C(L_{E_8}(k+2, 0),L_{E_8}(k, 0)\otimes L_{E_8}(1, 0)\otimes L_{E_8}(1, 0))$ is rational and $C_2$-cofinite.

On the other hand, $C(L_{E_8}(k+2, 0),L_{E_8}(k, 0)\otimes L_{E_8}(1, 0)\otimes L_{E_8}(1, 0))$ is also an extension of
$$C(L_{E_8}(k+2, 0),L_{E_8}(k, 0)\otimes L_{E_8}(2, 0))\bigotimes L(\frac{1}{2}, 0).$$
Moreover, by Lemma \ref{coset2}, $C(L_{E_8}(k+2, 0),L_{E_8}(k, 0)\otimes L_{E_8}(1, 0)\otimes L_{E_8}(1, 0))$ viewed as an $C(L_{E_8}(k+2, 0),L_{E_8}(k, 0)\otimes L_{E_8}(2, 0))\bigotimes L(\frac{1}{2}, 0)$-module has the following decomposition
$$C(L_{E_8}(k+2, 0),L_{E_8}(k, 0)\otimes L_{E_8}(2, 0))\otimes L(\frac{1}{2}, 0)\oplus W^1\otimes L(\frac{1}{2}, \frac{1}{2}) \oplus W^2\otimes L(\frac{1}{2}, \frac{1}{16}),$$
where $W^1, W^2$ denote some modules of $C(L_{E_8}(k+2, 0),L_{E_8}(k, 0)\otimes L_{E_8}(2, 0))$. Define a linear map $\sigma$ from $ C(L_{E_8}(k+2, 0),L_{E_8}(k, 0)\otimes L_{E_8}(1, 0)\otimes L_{E_8}(1, 0))$ to itself by
$$\sigma|_{W^2\otimes L(\frac{1}{2}, \frac{1}{16})}=-\id, ~~~\sigma|_{C(L_{E_8}(k+2, 0),L_{E_8}(k, 0)\otimes L_{E_8}(2, 0))\otimes L(\frac{1}{2}, 0)\oplus W^1\otimes L(\frac{1}{2}, \frac{1}{2})}=\id.$$
By Theorem 2.10 of \cite{ADL} and the fusion rules in Lemma \ref{fusion}, $\sigma$ is an automorphism of $ C(L_{E_8}(k+2, 0),L_{E_8}(k, 0)\otimes L_{E_8}(1, 0)\otimes L_{E_8}(1, 0))$. Hence, by Theorem 1 of \cite{M} and Theorem 5.24 of \cite{CM}, $C(L_{E_8}(k+2, 0),L_{E_8}(k, 0)\otimes L_{E_8}(2, 0))\otimes L(\frac{1}{2}, 0)\oplus W^1\otimes L(\frac{1}{2}, \frac{1}{2})$ is rational and $C_2$-cofinite. Define a linear map $\tau$ from $C(L_{E_8}(k+2, 0),L_{E_8}(k, 0)\otimes L_{E_8}(2, 0))\otimes L(\frac{1}{2}, 0)\oplus W^1\otimes L(\frac{1}{2}, \frac{1}{2})$ to itself by
$$\tau|_{C(L_{E_8}(k+2, 0),L_{E_8}(k, 0)\otimes L_{E_8}(2, 0))\otimes L(\frac{1}{2}, 0)}=\id, ~~~\tau|_{W^1\otimes L(\frac{1}{2}, \frac{1}{2})}=-\id.$$
Similarly, by Theorem 2.10 of \cite{ADL} and the fusion rules in Lemma \ref{fusion}, we can prove that $\tau$ is an automorphism of $C(L_{E_8}(k+2, 0),L_{E_8}(k, 0)\otimes L_{E_8}(2, 0))\otimes L(\frac{1}{2}, 0)\oplus W^1\otimes L(\frac{1}{2}, \frac{1}{2})$. Hence, by Theorem 1 of \cite{M} and Theorem 5.24 of \cite{CM}, $C(L_{E_8}(k+2, 0),L_{E_8}(k, 0)\otimes L_{E_8}(2, 0))\otimes L(\frac{1}{2}, 0)$ is rational and $C_2$-cofinite. Since $L(\frac{1}{2}, 0)$ is rational and $C_2$-cofinite, $C(L_{E_8}(k+2, 0),L_{E_8}(k, 0)\otimes L_{E_8}(2, 0))$ is also rational and $C_2$-cofinite (see \cite{DLM2} and \cite{DMZ}).
\qed

\section{Global dimensions and quantum dimensions of some diagonal coset vertex operator algebras}
\def\theequation{4.\arabic{equation}}
\setcounter{equation}{0}
Let $\g$ be a finite dimensional simple Lie algebra and $k,l$ be positive integers. In this section, we will calculate the global dimensions and quantum dimensions of the coset vertex operator algebra $C(L_{\g}(k+l,0),L_{\g}(k,0)\otimes L_{\g}(l,0))$ under the assumption that $C(L_{\g}(k+l,0),L_{\g}(k,0)\otimes L_{\g}(l,0))$ is rational and $C_2$-cofinite. As an application of the results in this section, we will classify irreducible modules of $C(L_{E_8}(k+2,0),L_{E_8}(k,0)\otimes L_{E_8}(2,0))$ in next section.

In the following, we will use $L_{\g}(k, \dot{\Lambda})$ and $L_{\g}(l, \ddot{\Lambda})$ to denote $L_{\g}(k, 0)$-module and $L_{\g}(l, 0)$-module, respectively. Note that $L_{\g}(k, \dot{\Lambda})\otimes L_{\g}(l, \ddot{\Lambda})$ may be viewed as an $L_{\g}(k+l, 0)$-module, and $L_{\g}(k+l, 0)$ is strongly regular. Hence, $L_{\g}(k, \dot{\Lambda})\otimes L_{\g}(l, \ddot{\Lambda})$ is completely reducible as an $L_{\g}(k+l, 0)$-module. For any $\Lambda\in P_+^{k+l}$, define $$M_{\dot{\Lambda}, \ddot{\Lambda}}^{\Lambda}=\Hom_{L_{\g}(k+l, 0)} (L_{\g}(k+l, \Lambda), L_{\g}(k, \dot{\Lambda})\otimes L_{\g}(l, \ddot{\Lambda})).$$
\begin{proposition}\cite{KW}
Let $k,l$ be positive integers. Then $M_{\dot{\Lambda}, \ddot{\Lambda}}^{\Lambda}\neq 0$ only if $\dot{\Lambda}+\ddot{\Lambda}-\Lambda\in Q$.
\end{proposition}

As a result, $L_{\g}(k, \dot{\Lambda})\otimes L_{\g}(l, \ddot{\Lambda})$ viewed as an $L_{\g}(k+l, 0)$-module has the following decomposition
$$L_{\g}(k, \dot{\Lambda})\otimes L_{\g}(l, \ddot{\Lambda})=\oplus_{\Lambda \in P_+^{k+l}; \dot{\Lambda}+\ddot{\Lambda}-\Lambda\in Q}L_{\g}(k+l, \Lambda)\otimes M_{\dot{\Lambda}, \ddot{\Lambda}}^{\Lambda}.$$
Note that  $$M_{\dot{\Lambda}, \ddot{\Lambda}}^{\Lambda}=\{v\in L_{\g}(k, \dot{\Lambda})\otimes L_{\g}(l, \ddot{\Lambda})| x\cdot v=0, h(0)\cdot v= \Lambda(h) v, \forall h\in \h, \forall x\in \hat{\g}_+\},$$ where $L_{\g}(k, \dot{\Lambda})\otimes L_{\g}(l, \ddot{\Lambda})$  is viewed as a module for the diagonal subalgebra of $\hat \g\oplus \hat \g$ and  $\hat{\g}_+=(\oplus_{\alpha\in \Delta_{\g}^+}\C x_\alpha)\bigoplus(\oplus_{x\in \g, n\in \Z_{>0}}\C x(n))$. Hence, $C(L_{\g}(k+l,0),L_{\g}(k,0)\otimes L_{\g}(l,0))=M_{0, 0}^0$ and $M_{\dot{\Lambda}, \ddot{\Lambda}}^{\Lambda}$ is a $C(L_{\g}(k+l,0),L_{\g}(k,0)\otimes L_{\g}(l,0))$-module.
\begin{proposition}\cite{KW}\label{positive}
Viewed as a $C(L_{\g}(k+l,0),L_{\g}(k,0)\otimes L_{\g}(l,0))$-module, the conformal weight of $M_{\dot{\Lambda}, \ddot{\Lambda}}^{\Lambda}$ is equal to or larger than $0$. Moreover, the conformal weight of $M_{\dot{\Lambda}, \ddot{\Lambda}}^{\Lambda}$ is equal to $0$ only if $(\dot{\Lambda}, \ddot{\Lambda},\Lambda)=(k\Lambda_i, l\Lambda_i,(k+l)\Lambda_i), i\in J$, where $\Lambda_0=0$.
\end{proposition}

We next show that there may be identifications between $C(L_{\g}(k+l,0),L_{\g}(k,0)\otimes L_{\g}(l,0))$-modules $\{M_{\dot{\Lambda}, \ddot{\Lambda}}^{\Lambda}|\dot\Lambda\in P_+^{k}, \ddot\Lambda\in P_+^{l}, \Lambda\in P_+^{k+l}\}$. First, let $h^i\in \h$ for $i=1, \cdots, l$ defined by $\alpha_i(h^j)=\delta_{i,j}$ for $j=1, \cdots, l$. For any $h\in \h$, set
$$\Delta(h, z)=z^{h(0)}\exp\left(\sum_{n=1}^{\infty}\frac{h(n)(-z)^{-n}}{-n}\right).$$
Given an $L_{\g}(k, 0)$-module $M$, we may construct a new $L_{\g}(k, 0)$-module $M^{(h^i)}$ by using the operator $\Delta(h^i, z)$. Explicitly, viewed as a vector space $M^{(h^i)}=M$, and the vertex operator map $Y_{M^{(h^i)}}(\cdot, z)=Y(\Delta(h^i, z)\cdot, z)$.
It was proved in \cite{L1} and \cite{L3} that for any $i\in J$ and any $\Lambda\in P_+^{k}$, $L_{\g}(k, \Lambda)^{(h^i)}$ is also an irreducible $L_{\g}(k, 0)$-module. Thus, there exists an element $\Lambda^{(i)}\in P_+^k$ such that  $L_{\g}(k, \Lambda)^{(h^i)}$ is isomorphic to $L_{\g}(k, \Lambda^{(i)})$.
\begin{theorem}\cite{L1,L3}
For $i\in J$, $L_{\g}(k, 0)^{(h^i)}\cong L_{\g}(k, k\Lambda_i)$.
\end{theorem}

As a consequence, we immediately obtain the following
\begin{corollary}\label{iden}
For any  $i\in J$ and $\dot\Lambda\in P_+^{k}, \ddot\Lambda\in P_+^{l}, \Lambda\in P_+^{k+l}$, we have $M_{\dot{\Lambda}, \ddot{\Lambda}}^{\Lambda}\cong M_{\dot{\Lambda}^{(i)}, \ddot{\Lambda}^{(i)}}^{\Lambda^{(i)}}$ as $C(L_{\g}(k+l,0),L_{\g}(k,0)\otimes L_{\g}(l,0))$-modules.
\end{corollary}

\begin{remark}
For $i\in J$ and $\dot\Lambda\in P_+^{k}, \ddot\Lambda\in P_+^{l}, \Lambda\in P_+^{k+l}$, the identification between $M_{\dot{\Lambda}, \ddot{\Lambda}}^{\Lambda}$ and $ M_{\dot{\Lambda}^{(i)}, \ddot{\Lambda}^{(i)}}^{\Lambda^{(i)}}$ was known in the physical reference \cite{SY}.
\end{remark}

To calculate the global dimension of $C(L_{\g}(k+l,0),L_{\g}(k,0)\otimes L_{\g}(l,0))$, we need the following results about $C(L_{\g}(k+l,0),L_{\g}(k,0)\otimes L_{\g}(l,0))$.
\begin{theorem}\label{crational}
Let $k, l$ be positive integers. Suppose that the vertex operator algebra $C(L_{\g}(k+l,0),L_{\g}(k,0)\otimes L_{\g}(l,0))$ is rational and $C_2$-cofinite. Then\\
(1) $C(L_{\g}(k+l,0),L_{\g}(k,0)\otimes L_{\g}(l,0))$ is strongly regular. \\
(2) Any irreducible $C(L_{\g}(k+l,0),L_{\g}(k,0)\otimes L_{\g}(l,0))$-module is isomorphic to a submodule of $M_{\dot{\Lambda}, \ddot{\Lambda}}^{\Lambda}$ for some $\dot\Lambda\in P_+^{k}, \ddot\Lambda\in P_+^{l}, \Lambda\in P_+^{k+l}$.\\
(3) All the conformal weights of irreducible $C(L_{\g}(k+l,0),L_{\g}(k,0)\otimes L_{\g}(l,0))$-modules except $C(L_{\g}(k+l,0),L_{\g}(k,0)\otimes L_{\g}(l,0))$ are larger than $0$.
\end{theorem}
\pf (1) By assumption, we only need to $C(L_{\g}(k+l,0),L_{\g}(k,0)\otimes L_{\g}(l,0))$ is self-dual. By Corollary 3.2 of \cite{Li} and Lemma \ref{scoset}, it is good enough to prove that $(\w^1+\w^2-\w^3)_2v=0$ for any $v\in C(L_{\g}(k+l,0),L_{\g}(k,0)\otimes L_{\g}(l,0))_1$. By the definition of $C(L_{\g}(k+l,0),L_{\g}(k,0)\otimes L_{\g}(l,0))$, we have $\w^3_2v=0$. Moreover,  $\w^3_1v=0$, this implies that $v$ is an element of $(L_{\g}(k,0)\otimes L_{\g}(l,0))_1$. Hence,  by the formula (6.2.45) of \cite{LL}, $(\w^1+\w^2)_2v =0$. Therefore, $(\w^1+\w^2-\w^3)_2v=0$, as desired.

(2) This follows immediately from Theorem 2 of \cite{KM}.

(3) This follows immediately from Proposition \ref{positive} and Corollary \ref{iden}.
\qed

\vskip.25cm
We are now ready to prove the main result in this section.
\begin{theorem}\label{qdimc}
Let $k, l$ be positive integers. Suppose that the vertex operator algebra $C(L_{\g}(k+l,0),L_{\g}(k,0)\otimes L_{\g}(l,0))$ is rational and $C_2$-cofinite. Then we have\\
(1) $\Glob C(L_{\g}(k+l,0),L_{\g}(k,0)\otimes L_{\g}(l,0))=\frac{1}{|P/Q|^2S_{L_{\g}(k, 0), L_{\g}(k, 0)}^2S_{L_{\g}(l, 0), L_{\g}(l, 0)}^2S_{L_{\g}(k+l, 0), L_{\g}(k+l, 0)}^2 }$.\\
(2) For any $\dot\Lambda\in P_+^{k}, \ddot\Lambda\in P_+^{l}, \Lambda\in P_+^{k+l}$ such that $\dot{\Lambda}+\ddot{\Lambda}-\Lambda\in Q$, we have $$\qdim_{C(L_{\g}(k+l,0),L_{\g}(k,0)\otimes L_{\g}(l,0))} M_{\dot{\Lambda}, \ddot{\Lambda}}^{\Lambda}=\frac{S_{L_{\g}(k, 0), L_{\g}(k, \dot\Lambda)}S_{L_{\g}(l, 0), L_{\g}(l, \ddot\Lambda)}S_{L_{\g}(k+l, 0), L_{\g}(k+l, \Lambda)}}{S_{L_{\g}(k, 0), L_{\g}(k, 0)}S_{L_{\g}(l, 0), L_{\g}(l, 0)}S_{L_{\g}(k+l, 0), L_{\g}(k+l, 0)} }.$$
\end{theorem}
\pf (1) By the assumption and Theorem \ref{crational}, all the conditions in Theorem \ref{qdim1} hold for $C(L_{\g}(k+l,0),L_{\g}(k,0)\otimes L_{\g}(l,0))$. Let $C(L_{\g}(k+l,0),L_{\g}(k,0)\otimes L_{\g}(l,0))=M^0, M^1, \cdots, M^p$ be all the non-isomorphic irreducible $C(L_{\g}(k+l,0),L_{\g}(k,0)\otimes L_{\g}(l,0))$-modules. Then there exist complex numbers $S_{i, j}, i,j=0,\cdots, p$, such that
$$Z_{M^i}(v, \frac{-1}{\tau})=\tau^{\wt [v]}\sum_{j=0}^p S_{i,j} Z_{M^j}(v, \tau).$$
Moreover, the matrix $(S_{i, j})$ is independent of the choice of $v\in C(L_{\g}(k+l,0),L_{\g}(k,0)\otimes L_{\g}(l,0))$. In particular, we have
\begin{equation}\label{S1}
Z_{M^i}(\1, \frac{-1}{\tau})=\sum_{j=0}^p S_{i,j} Z_{M^j}(\1, \tau).
\end{equation}
On the other hand, by the formulas (12.2.2) of \cite{K} and (2.7.2) of \cite{KW}, we have
\begin{align*}
&Z_{M_{L_{\g}(k, 0), L_{\g}(l, 0)}^{L_{\g}(k+l, 0)}}(\1, \frac{-1}{\tau})\\
&=\sum_{\substack{\dot\Lambda\in P_+^{k}, \ddot\Lambda\in P_+^{l}, \Lambda\in P_+^{k+l},\\ \dot{\Lambda}+\ddot{\Lambda}-\Lambda\in Q}} S_{L_{\g}(k, 0), L_{\g}(k, \dot\Lambda)}S_{L_{\g}(l, 0), L_{\g}(l, \ddot\Lambda)}S_{L_{\g}(k+l, 0), L_{\g}(k+l, \Lambda)} Z_{M_{L_{\g}(k, \dot\Lambda), L_{\g}(l, \ddot\Lambda)}^{L_{\g}(k+l, \Lambda)}}(\1, \tau).
\end{align*}
Comparing with the formula (\ref{S1}), it follows from  Proposition \ref{positive} and Theorem \ref{crational} that
$$S_{0, 0}=|P/ Q|S_{L_{\g}(k, 0), L_{\g}(k, 0)}S_{L_{\g}(l, 0), L_{\g}(l, 0)}S_{L_{\g}(k+l, 0), L_{\g}(k+l, 0)}.$$
Furthermore, by Theorems \ref{qdim1}, \ref{crational}, we have
\begin{align*}
&\Glob C(L_{\g}(k+l,0),L_{\g}(k,0)\otimes L_{\g}(l,0))\\
&=\frac{1}{S_{0, 0}^2}=\frac{1}{|P/Q|^2S_{L_{\g}(k, 0), L_{\g}(k, 0)}^2S_{L_{\g}(l, 0), L_{\g}(l, 0)}^2S_{L_{\g}(k+l, 0), L_{\g}(k+l, 0)}^2 }.
\end{align*}
(2) By the definition of  quantum dimension, we have
\begin{align*}
\qdim_{C(L_{\g}(k+l,0),L_{\g}(k,0)\otimes L_{\g}(l,0))} M_{\dot{\Lambda}, \ddot{\Lambda}}^{\Lambda}=\lim_{y\to 0^+}\frac{Z_{M_{\dot{\Lambda}, \ddot{\Lambda}}^{\Lambda}}(\1, \sqrt{-1}y)}{Z_{C(L_{\g}(k+l,0),L_{\g}(k,0)\otimes L_{\g}(l,0))}(\1, \sqrt{-1}y)},
\end{align*}
where $y$ is real and positive. On the other hand, for any $\dot\Lambda\in P_+^{k}, \ddot\Lambda\in P_+^{l}, \Lambda\in P_+^{k+l}$ such that $\dot{\Lambda}+\ddot{\Lambda}-\Lambda\in Q$, it is known by the formula (2.7.15) of \cite{KW} that
$$\lim_{y\to 0^+} \frac{Z_{M_{\dot{\Lambda}, \ddot{\Lambda}}^{\Lambda}}(\1, \sqrt{-1}y)}{|P/Q|S_{L_{\g}(k, 0), L_{\g}(k, \dot\Lambda)}S_{L_{\g}(l, 0), L_{\g}(l, \ddot\Lambda)}S_{L_{\g}(k+l, 0), L_{\g}(k+l, \Lambda)}e^{\pi (z_k+z_l-z_{k+l})/12y}}=1,$$ where $z_{n}=\frac{(\dim \g)n}{n+h^{\vee}}$ for $n\in \Z_{>0}$. It follows that
\begin{align*}
&\qdim_{C(L_{\g}(k+l,0),L_{\g}(k,0)\otimes L_{\g}(l,0))} M_{\dot{\Lambda}, \ddot{\Lambda}}^{\Lambda}\\
&=\lim_{y\to 0^+}\frac{|P/Q|S_{L_{\g}(k, 0), L_{\g}(k, \dot\Lambda)}S_{L_{\g}(l, 0), L_{\g}(l, \ddot\Lambda)}S_{L_{\g}(k+l, 0), L_{\g}(k+l, \Lambda)}e^{\pi (z_k+z_l-z_{k+l})/12y}}{|P/Q|S_{L_{\g}(k, 0), L_{\g}(k, 0)}S_{L_{\g}(l, 0), L_{\g}(l, 0)}S_{L_{\g}(k+l, 0), L_{\g}(k+l, 0)}e^{\pi (z_k+z_l-z_{k+l})/12y}}\\
&=\frac{S_{L_{\g}(k, 0), L_{\g}(k, \dot\Lambda)}S_{L_{\g}(l, 0), L_{\g}(l, \ddot\Lambda)}S_{L_{\g}(k+l, 0), L_{\g}(k+l, \Lambda)}}{S_{L_{\g}(k, 0), L_{\g}(k, 0)}S_{L_{\g}(l, 0), L_{\g}(l, 0)}S_{L_{\g}(k+l, 0), L_{\g}(k+l, 0)}},
\end{align*}
as desired.
\qed
\section{Irreducible modules of $C(L_{E_8}(k+2,0), L_{E_8}(k,0)\otimes L_{E_8}(2,0))$}
As an application of the results obtained in Section 4, we classify irreducible modules of the coset vertex operator algebra $C(L_{E_8}(k+2,0), L_{E_8}(k,0)\otimes L_{E_8}(2,0))$ in this section. First, we prove the following general result.
\begin{theorem}\label{general}
Let $\g$ be a finite dimensional simple Lie algebra and  $k, l$ be positive integers. Suppose the following two conditions hold:\\
(i) The vertex operator algebra $C(L_{\g}(k+l,0),L_{\g}(k,0)\otimes L_{\g}(l,0))$ is rational and $C_2$-cofinite.\\
(ii)  There exist $C(L_{\g}(k+l,0),L_{\g}(k,0)\otimes L_{\g}(l,0))$-modules $W^1, \cdots, W^s$ such that for any $\dot\Lambda\in P_+^{k}, \ddot\Lambda\in P_+^{l}, \Lambda\in P_+^{k+l}$, $ M_{\dot{\Lambda}, \ddot{\Lambda}}^{\Lambda}$ is a direct sum of some $W^{i_1}, \cdots, W^{i_t}$. Moreover, $$\sum_{i=1}^s (\qdim_{C(L_{\g}(k+l,0),L_{\g}(k,0)\otimes L_{\g}(l,0))} W^i)^2=\Glob C(L_{\g}(k+l,0),L_{\g}(k,0)\otimes L_{\g}(l,0)).$$
Then $W^1, \cdots, W^s$ are irreducible $C(L_{\g}(k+l,0),L_{\g}(k,0)\otimes L_{\g}(l,0))$-modules. Moreover, $W^1, \cdots, W^s$ are all the non-isomorphic irreducible $C(L_{\g}(k+l,0),L_{\g}(k,0)\otimes L_{\g}(l,0))$-modules.
\end{theorem}
\pf By the condition (i) and Theorem \ref{crational}, all the conditions in Theorem \ref{qdim1} hold for $C(L_{\g}(k+l,0),L_{\g}(k,0)\otimes L_{\g}(l,0))$. Let $C(L_{\g}(k+l,0),L_{\g}(k,0)\otimes L_{\g}(l,0))=M^0, M^1, \cdots, M^p$ be all the non-isomorphic irreducible $C(L_{\g}(k+l,0),L_{\g}(k,0)\otimes L_{\g}(l,0))$-modules. By Theorem \ref{crational}, each $M^i$  is isomorphic to a submodule of $ M_{\dot{\Lambda}, \ddot{\Lambda}}^{\Lambda}$ for some $\dot\Lambda\in P_+^{k}, \ddot\Lambda\in P_+^{l}, \Lambda\in P_+^{k+l}$. It follows from  the condition (ii) that each $M^i$ is isomorphic to a submodule of $W^j$ for some $j$. By Theorem \ref{qdim1}, $\qdim_{C(L_{\g}(k+l,0),L_{\g}(k,0)\otimes L_{\g}(l,0))} M^i\geq 1 $. On the other hand, by the condition (ii), we have
\begin{align*}
\sum_{i=1}^s (\qdim_{C(L_{\g}(k+l,0),L_{\g}(k,0)\otimes L_{\g}(l,0))} W^i)^2&=\Glob C(L_{\g}(k+l,0),L_{\g}(k,0)\otimes L_{\g}(l,0))\\
&=\sum_{i=0}^{p}(\qdim_{C(L_{\g}(k+l,0),L_{\g}(k,0)\otimes L_{\g}(l,0))} M^i)^2.
\end{align*}
It follows that each $W^j$ must be irreducible. Moreover, $W^1, \cdots, W^s$ must be non-isomorphic. Thus, $W^1, \cdots, W^s$ are all the non-isomorphic $C(L_{\g}(k+l,0),L_{\g}(k,0)\otimes L_{\g}(l,0))$-modules.
\qed

\vskip.25cm
Finally, we classify irreducible modules of the coset vertex operator algebra $C(L_{E_8}(k+2,0), L_{E_8}(k,0)\otimes L_{E_8}(2,0))$. We will need the following fact which was obtained in Corollary 2.7 of \cite{KW}.
\begin{proposition}\label{S3}
Let $\g$ be a finite dimensional simple Lie algebra and $k$ be a positive integer. Then the sum of $S_{L_{\g}(k, 0), L_{\g}(k, \Lambda)}^2$, where $\Lambda$ runs over a congruence class of $P_+^k$ mod $Q$, is equal to $|P/Q|^{-1}$.
\end{proposition}
We are now ready to classify irreducible modules of $C(L_{E_8}(k+2,0), L_{E_8}(k,0)\otimes L_{E_8}(2,0))$.
 \begin{theorem}\label{main2}
 Let $k$ be a positive integer. Then  $$\{ M_{\dot{\Lambda}, \ddot{\Lambda}}^{\Lambda}|\dot\Lambda\in P_+^{k}, \ddot\Lambda\in P_+^{2}, \Lambda\in P_+^{k+2}, \dot{\Lambda}+\ddot{\Lambda}-\Lambda\in Q\}$$
 are all the non-isomorphic irreducible $C(L_{E_8}(k+2,0), L_{E_8}(k,0)\otimes L_{E_8}(2,0))$-modules.
 \end{theorem}
 \pf By Theorem \ref{main1}, the vertex operator algebra $C(L_{E_8}(k+2,0), L_{E_8}(k,0)\otimes L_{E_8}(2,0))$ is rational and $C_2$-cofinite.
 We next verify the condition (ii) of Theorem \ref{general}. Recall that when $\g=E_8$,  we have $|P/Q|=1$ \cite{KW}. Moreover, by the formula (2.2.4) of \cite{KW}, for any positive integer $n$,
 $$\sum_{\Lambda\in P_+^{n}}S_{L_{\g}(n, 0), L_{\g}(n, \Lambda)}^2=1.$$
 Hence, by Theorem \ref{qdimc} and Proposition \ref{S3}, we have
 \begin{align*}
 &\sum_{\substack{\dot\Lambda\in P_+^{k}, \ddot\Lambda\in P_+^{2}, \Lambda\in P_+^{k+2},\\\dot{\Lambda}+\ddot{\Lambda}-\Lambda\in Q}}(\qdim_{C(L_{E_8}(k+2,0), L_{E_8}(k,0)\otimes L_{E_8}(2,0))} M_{\dot{\Lambda}, \ddot{\Lambda}}^{\Lambda})^2\\
 &=\sum_{\substack{\dot\Lambda\in P_+^{k}, \ddot\Lambda\in P_+^{2}, \Lambda\in P_+^{k+2},\\\dot{\Lambda}+\ddot{\Lambda}-\Lambda\in Q}}\left(\frac{S_{L_{E_8}(k, 0), L_{E_8}(k, \dot\Lambda)}S_{L_{E_8}(2, 0), L_{E_8}(2, \ddot\Lambda)}S_{L_{E_8}(k+2, 0), L_{E_8}(k+2, \Lambda)}}{S_{L_{E_8}(k, 0), L_{E_8}(k, 0)}S_{L_{E_8}(2, 0), L_{E_8}(2, 0)}S_{L_{E_8}(k+2, 0), L_{E_8}(k+2, 0)} }\right)^2\\
 &=\sum_{\dot\Lambda\in P_+^{k}, \ddot\Lambda\in P_+^{2}}\left(\frac{S_{L_{E_8}(k, 0), L_{E_8}(k, \dot\Lambda)}S_{L_{E_8}(2, 0), L_{E_8}(2, \ddot\Lambda)}}{S_{L_{E_8}(k, 0), L_{E_8}(k, 0)}S_{L_{E_8}(2, 0), L_{E_8}(2, 0)}S_{L_{E_8}(k+2, 0), L_{E_8}(k+2, 0)} }\right)^2\\
 &=\left(\frac{1}{S_{L_{E_8}(k, 0), L_{E_8}(k, 0)}S_{L_{E_8}(2, 0), L_{E_8}(2, 0)}S_{L_{E_8}(k+2, 0), L_{E_8}(k+2, 0)} }\right)^2\\
 &=\Glob C(L_{E_8}(k+2,0), L_{E_8}(k,0)\otimes L_{E_8}(2,0)).
 \end{align*}
 Therefore, by Theorem \ref{general}, $$\{ M_{\dot{\Lambda}, \ddot{\Lambda}}^{\Lambda}|\dot\Lambda\in P_+^{k}, \ddot\Lambda\in P_+^{2}, \Lambda\in P_+^{k+2}, \dot{\Lambda}+\ddot{\Lambda}-\Lambda\in Q\}$$
 are all the non-isomorphic irreducible $C(L_{E_8}(k+2,0), L_{E_8}(k,0)\otimes L_{E_8}(2,0))$-modules.
 \qed

\end{document}